\documentclass[12pt]{amsart}

\usepackage{vmargin}
\usepackage[dvips]{graphicx}
\usepackage{color}
\usepackage{xypic}
\xyoption{all}
\input{amssym.def}
\input{amssym}
\usepackage{a4wide}
\usepackage{supertabular}
\usepackage{array}
\usepackage{multicol}

\setmarginsrb{3cm}{2cm}{3cm}{3cm}{75pt}{20pt}{20pt}{30mm}
\newcommand{\A}{\mathbb A}

 \newcommand{\G}{\mathbb{G}}
 \newcommand{\Gm}{\mathbb{G}_{\mathrm {m}}}

 \newcommand{\PP}{\mathbb{P}}
 \renewcommand{\P}{\PP}

 \newcommand{\GL}{\mathrm{GL}}

  \newcommand{\End}{\mathrm{End}}

 \newcommand{\Hom}{\mathrm{Hom}}
 
\newcommand{\ddim}{\mathrm{dim}}

 \newcommand{\PGL}{\mathrm{PGL}}

 \theoremstyle{plain}
 \newtheorem{thm}{Th\'eor\`eme}[section]
 \newtheorem{defi}[thm]{D\'efinition}

 \newtheorem{prop}[thm]{Proposition}
 \newtheorem{lem}[thm]{Lemme}
 \newtheorem{coro}[thm]{Corollaire}

 \theoremstyle{remark}
 \newtheorem{rem}[thm]{Remarque}

 \newenvironment{dem}{{\bf D\'emonstration.}}{\hfill$\square$}

\def\signmf{\bigskip \begin{center} {\sc Mathieu Florence\par\vspace{3mm}
 Universit\'e Pierre et Marie Curie \par Institut de Math\'ematiques de Jussieu, Equipe de Topologie et G\'eom\'etrie Alg\'ebriques \par
175, rue du Chevaleret\par
75013 Paris \par\vspace{3mm}
e-mail:} \tt{florence@math.jussieu.fr} \end{center}}

\begin{document}

\title{G\'eom\'etrie birationnelle \'equivariante des grassmanniennes}

\author{M. Florence}

\maketitle

\vspace*{-10mm}

\begin{abstract} 
Soit $K$ un corps infini et $A$ une $K$-alg\`ebre de dimension finie $n$. Soit $0<r<n$ un entier. Le  r\'esultat principal de cet article est le suivant. Sous certaines hypoth\`eses sur $A$ (satisfaites si $A$ est \'etale), la grassmannienne $\G(r,A)$ est $K$-birationnelle, de mani\`ere $\PGL_1(A)$-\'equivariante, au produit de  $\G(\mathrm{pgcd}(r,n),A)$ par un espace affine (th\'eor\`eme \ref{thprinc}). On en d\'eduit par torsion plusieurs r\'esultats nouveaux sur les vari\'et\'es de Severi-Brauer g\'en\'eralis\'ees, li\'es \`a la conjecture d'Amitsur. Les plus significatifs sont les suivants.\\
\noindent i) Si $A$ est une $K$-alg\`ebre simple centrale de degr\'e $n$, et si $0<r<n$ est un entier, la $r$-i\`eme vari\'et\'e de Severi-Brauer g\'en\'eralis\'ee $SB(r,A)$ est $K$-birationnelle au produit de $SB(\mathrm{pgcd}(r,n),A)$ par un $K$-espace affine de dimension convenable.\\
\noindent ii) Si $A$ et $B$ sont deux $K$-alg\`ebres simples centrales de degr\'es premiers entre eux, $SB(A \otimes_K B)$ est birationnelle au produit de $SB(A) \times_K SB(B)$ par un $K$-espace affine de dimension convenable. \\
\end{abstract}

%\tableofcontents

{\bf Mots-cl\'es.} G\'eom\'etrie birationnelle, grassmanniennes, vari\'et\'es de Severi-Brauer, conjecture d'Amitsur.

\newpage
\section*{Introduction}
\bigbreak%%%
Soit $K$ un corps. En 1955, Amitsur a formul\'e la conjecture suivante  (\cite{Ami}) . Soient $A$, $B$ deux $K$-alg\`ebres simples centrales de m\^eme degr\'e $n$ (racine carr\'ee de la dimension sur $K$). Alors les vari\'et\'es de Severi-Brauer associ\'ees \`a $A$ et $B$ (not\'ees respectivement $SB(A)$ et $SB(B)$) sont $K$-birationnelles si et seulement si les sous-groupes cycliques engendr\'es par les classes de $A$ et de $B$ dans le groupe de Brauer de $K$ co\"incident. Amitsur a d\'emontr\'e l'implication `seulement si' de sa conjecture. En d\'epit de r\'esultats partiels obtenus notamment par Roquette, Tregub et Krashen (cf. par exemple \cite{Kra} pour une exposition de ces r\'esultats), cette conjecture reste encore largement ouverte. Il n'est cependant pas difficile de d\'emontrer que, si les sous-groupes engendr\'es par les classes de $A$ et de $B$ dans le groupe de Brauer de $K$ co\"incident, alors $SB(A)$ et $SB(B)$ sont \textit{stablement} birationnellement isomorphes. En effet, soit $E$ le corps des fonctions de $SB(A)$. L'extension $E/K$ d\'eploie $A$, donc aussi $B$. Par suite, $SB(B)\times_K E$ est $E$-isomorphe \`a l'espace projectif $\P^{n-1}_E$. Donc $SB(B) \times_K SB(A)$ est $K$-birationnellement isomorphe \`a $\P^{n-1}_K \times_K SB(A)$. En inversant les r\^oles de $A$ et de $B$, on obtient que $\P^{n-1}_K \times_K SB(A)$ est $K$-birationnelle \`a $\P^{n-1}_K \times_K SB(B)$, cqfd. En g\'en\'eral, nous dirons ici que deux $K$-vari\'et\'es $X$ et $Y$ sont stablement $K$-birationnelles s'il existe deux entiers $n$ et $m$ tels que $X \times_K \P_K^n$ est $K$-birationnelle \`a $Y \times_K \P_K^m$. En s'inspirant de la preuve pr\'ec\'edente, le lecteur pourra d\'emontrer les propositions suivantes.\\

i) Soient $A$, $B$ deux $K$-alg\`ebres simples centrales de degr\'es premiers entre eux. Alors $SB(A\otimes_K B)$ est stablement $K$-birationnelle \`a $SB(A) \times_K SB(B)$.\\

ii) Soit $A$ une $K$-alg\`ebre simple centrale de degr\'e $n$. Soit $r < n$ un entier positif. Soit $SB(r,A)$ la $r$-i\`eme vari\'et\'e de Severi-Brauer g\'en\'eralis\'ee associ\'ee \`a $A$; c'est la $K$-forme de la grassmannienne $\G(r,n)$ qui repr\'esente les id\'eaux \`a droite de $A$, de dimension $rn$. Alors $SB(r,A)$ et $SB(pgcd(n,r),A)$ sont stablement $K$-birationnelles.\\

\noindent On peut formuler les propri\'et\'es plus fines suivantes, et se demander si elles sont vraies. \\

i') Soient $A$, $B$ deux $K$-alg\`ebres simples centrales de degr\'es premiers entre eux. Alors $SB(A\otimes_K B)$ est $K$-birationnelle au produit de $SB(A) \times_K SB(B)$ par un espace projectif (ou affine, c'est pareil) de dimension convenable.\\

ii') M\^emes hypoth\`eses que ii), avec la conclusion: $SB(r,A)$ est $K$-birationnelle au produit de $SB(pgcd(n,r),A)$ par un espace affine de dimension convenable.\\

De m\^eme que la conjecture d'Amitsur, le passage de i) \`a i') (resp. de ii) \`a ii')) peut \^etre vu comme un cas particulier d'un probl\`eme de `simplification par l'espace affine': si $X$ et $Y$ sont deux $K$-vari\'et\'es telles que $X \times_K \A_K^1$ est $K$-birationnelle \`a $Y \times_K \A_K^1$, est-il vrai que $X$ et $Y$ sont $K$-birationnelles? M\^eme si $X$ est $K$-rationnelle, on sait que la r\'eponse \`a cette question, qui n'est alors autre que la conjecture de Zariski, est en g\'en\'eral n\'egative (\cite{BCSS}).\\
Dans cette article, nous \'etudions la g\'eom\'etrie birationnelle des grassmanniennes de $K$-alg\`ebres de dimension finie. Plus pr\'ecis\'ement, soit $A$ une $K$-alg\`ebre de dimension finie. Pour tout entier $r$ compris entre $0$ et la dimension de $A$, la grassmannienne $\G(r,A)$ des $r$-plans de $A$ est naturellement munie d'une action du $K$-groupe alg\'ebrique $\GL_1(A)$, dont les $K$-points sont les \'el\'ements inversibles de $A$. On s'int\'eresse, par une approche nouvelle \`a la connaissance de l'auteur, \`a certaines propri\'et\'es birationnelles de cette action. Le r\'esultat principal est le th\'eor\`eme \ref{thprinc}. Il affirme que, sous certaines hypoth\`eses sur $A$  (v\'erifi\'ees si $A/K$ est \'etale), la $K$-vari\'et\'e  $\G(r,A)$ est $K$-birationnelle, de mani\`ere $\GL_1(A)$-\'equivariante, au produit de $\G(pgcd(r,n),A)$ par un espace affine sur lequel $\GL_1(A)$ agit trivialement. Les isomorphismes birationnels obtenus sont non triviaux. Leur complexit\'e d\'epend de fa\c{c}on cruciale du nombre d'\'etapes de l'algorithme d'Euclide, prenant en entr\'ee les entiers $n$ et $r$, et donnant leur plus grand commun diviseur en sortie.\\
\noindent L'article est structur\'e comme suit. Dans une premi\`ere partie, nous introduisons la notion de $\textit{bonne}$ $K$-alg\`ebre (d\'efinition \ref{bonne}) ainsi que les vari\'et\'es auxiliaires $G(r,s,A)$ associ\'ees \`a une bonne $K$-alg\`ebre $A$ (d\'efinition \ref{defG}). On d\'egage ensuite certaines relations birationnelles et \'equivariantes entre ces vari\'et\'es (d\'efinition \ref{defiphi}), qui motivent a posteriori leur introduction. Dans une seconde partie, nous exploitons ces relations pour en d\'eduire le th\'eor\`eme principal (th\'eor\`eme \ref{thprinc}). Dans la derni\`ere section, nous exposons quelques applications, en majeure partie nouvelles, de ce th\'eor\`eme. Entre autres, nous d\'emontrons  la validit\'e de i') et de ii'), ainsi que, dans certains cas, la conjecture d'Amitsur g\'en\'eralis\'ee \'enonc\'ee par Krashen dans \cite{Kra}.

\section{Notations, rappels}
Dans toute la suite, on d\'esigne par $K$ un corps infini, et par $\overline K$ une cl\^oture alg\'ebrique de $K$. Introduisons un bref dictionnaire destin\'e \`a all\'eger les \'enonc\'es. Ainsi, les termes `vectoriel', `alg\`ebre', `vari\'et\'e', `vari\'et\'e rationnelle', `groupe alg\'ebrique' signifient respectivement, sauf mention du contraire, `$K$-espace vectoriel de dimension finie', `$K$-alg\`ebre de dimension finie', `$K$-sch\'ema de type fini', `vari\'et\'e $K$-rationnelle', `$K$-groupe alg\'ebrique lin\'eaire'. Si $X$ et $Y$ sont deux vari\'et\'es (resp. vectoriels), on note par $X \times Y$ (resp. $X \otimes Y$) la vari\'et\'e $X \times_K Y$ (resp. le vectoriel $X \otimes_K Y$).
Soit $V$ un vectoriel. On note $V^*$ le dual de $V$, et $\A (V)$ l'espace affine de $V$; c'est la vari\'et\'e dont le foncteur des points est $\A (V)(A)=V\otimes_K A$ pour toute alg\`ebre commutative $A$ (non n\'ecessairement de dimension finie). Si $m$ est un entier positif, on note $\G(m,V)$ la grassmannienne des $m$-sous-vectoriels de $V$. On rappelle que sa dimension est $m(\ddim V -m)$. L'espace projectif de $V$ est $\G(1,V)$; on le d\'esigne par $\P(V)$. Si $E \subset V$ est un sous-vectoriel, on pose $E^\perp= \{ \phi \in V^*; \phi(E)=0\}$. Si $A$ est une alg\`ebre, on note $\GL_1(A)$ le groupe alg\'ebrique dont l'ensemble des points dans une alg\`ebre commutative $B$ est le groupe des \'el\'ements inversibles de $A \otimes_K B$. On note $\PGL_1(A)$ le quotient $ \GL_1(A)/ \Gm$; c'est l'ouvert de $\P(A)$ form\'e des droites dont tout \'el\'ement non nul est inversible pour la multiplication de $A$. Le groupe alg\'ebrique $\GL_1(A)$ poss\`ede deux repr\'esentations \'evidentes qui seront consid\'er\'ees dans la suite. D'une part, il agit par multiplication \`a gauche sur le vectoriel $A$. D'autre part, il agit sur le vectoriel $A^*$ par la formule $$(a.\phi)(x)=\phi(a^{-1}x),$$ pour $a \in \GL_1(A)(K)$, $\phi \in A^*$ et $x \in A$. Pour tout entier $m$, ces repr\'esentations induisent une action naturelle de $\PGL_1(A)$ sur  $\G(m,A)$ et sur $\G(m,A^*)$, capitale dans ce qui suit. Enfin, nous noterons les applications rationnelles par des fl\`eches pleines, et non par des pointill\'es, comme l'on a coutume de faire. Ceci est justifi\'e par le fait que la quasi-totalit\'e des fl\`eches consid\'er\'ees dans cet article sont des applications rationnelles, et ne saurait pr\^eter \`a confusion.

\section{Quelques outils}

\begin{defi}\label{bonouvert}
Soit $G$ un groupe alg\'ebrique, agissant sur une vari\'et\'e non vide $X$ de fa\c{c}on g\'en\'eriquement libre (i.e. il existe un ouvert Zariski dense de $X$, $G$-stable, o\`u l'action est libre). Par un th\'eor\`eme de Gabriel (\cite{SGA3}, Expos\'e V, Th\'eor\`eme 8.1), il existe alors un ouvert Zariski dense $U \subset X$, $G$-stable, tel que le quotient $U \longrightarrow U/G$ existe dans la cat\'egorie des vari\'et\'es, et est un $G$-torseur. Un ouvert $U$ poss\'edant cette propri\'et\'e sera appel\'e bon pour l'action de $G$ sur $X$. 
\end{defi}
\begin{defi}
Soit $A$ une alg\`ebre et $M$ un $A$-module \`a gauche. Soit $E$ (resp. $X$) un sous-vectoriel de $A$ (resp. de $M$). Le produit de $X$ par $E$, not\'e $E.X$, est par d\'efinition l'image de la compos\'ee $$E \otimes X \longrightarrow A \otimes M \longrightarrow M,$$ o\`u le premier morphisme est le produit tensoriel des deux inclusions, et le second la loi donnant l'action de $A$ sur $M$. On d\'efinit de mani\`ere semblable un produit $X.E$ pour les $A$-modules \`a droite.
\end{defi}
Nous introduisons maintenant les vari\'et\'es auxiliaires $G'(r,s,U)$. Ces vari\'et\'es jouent un r\^ole crucial dans la strat\'egie que nous adpotons pour d\'emontrer le th\'eor\`eme \ref{thprinc}. Leur int\'er\^et peut sembler obscur \`a ce stade. Il deviendra, nous l'esp\'erons, plus clair au fil de la lecture de cet article. 
\begin{defi} \label{defG'}
Soit $A$ une alg\`ebre de dimension $n$. Notons par $$<.,.> : A \times A^* \longrightarrow K$$ l'accouplement canonique. Soient $r,s,u$ trois entiers positifs tels que $n \geq su+r$ et $n \geq ru+s$. Soit $U \subset A$ un  sous-vectoriel de dimension $u$. On appelle $G'(r,s,U)$ la sous-vari\'et\'e ferm\'ee de $\G(r,A^*) \times \G(s,A)$ donn\'ee par $$\{ (X,Y) \in \G(r,A^*) \times \G(s,A); <Y.U,X>=0 \}.$$ 
\end{defi}

\begin{lem} \label{lemG'}
On conserve les objets et les notations de la d\'efinition \ref{defG'}. Les assertions suivantes sont vraies.\\

\noindent i) $G'(r,s,U)$ est munie d'une action naturelle de $\PGL_1(A)$. \\
\noindent ii) La fibre de la premi\`ere projection $ G'(r,s,U) \longrightarrow \G(r,A^*)$ (resp. de la seconde projection $ G'(r,s,U) \longrightarrow \G(s,A)$) en $X \in \G(r,A^*)(K)$ (resp. $Y \in \G(s,A)(K)$) est $\G(s, (U.X) ^\perp)$ (resp. $\G(r, (Y.U) ^\perp)$). On a l'assertion analogue en rempla\c{c}ant $K$ par une extension quelconque de $K$.\\
\noindent iii) Les deux projections pr\'ec\'edentes sont dominantes.
\end{lem}
\begin{dem}
Pour le premier point, on constate tout de suite que l'action naturelle de $\PGL_1(A)$ sur $\G(r,A^*) \times \G(s,A)$ laisse $ G'(r,s,U)$ stable. L'assertion ii) est imm\'ediate. On a $\dim \G(s, (U.X) ^\perp)= s(\dim (U.X) ^\perp -s) \geq s(n-ru-s) \geq 0$, et l'in\'egalit\'e analogue pour $\G(r, (Y.U) ^\perp)$, d'o\`u iii). 
\end{dem}

Nous pouvons maintenant d\'efinir la notion de bonne alg\`ebre, capitale dans ce qui suit.

\begin{defi} \label{bonne}
Soit $A$ une alg\`ebre, dont on note $n$ la dimension. On dit qu'un sous-vectoriel $U$ de $A$, de dimension $u$, est bon si la propri\'et\'e qui suit est vraie. Pour tous entiers positifs $r$ et $s$ tels que $n\geq ru+s$ et $n\geq su+r$, il existe $(X,Y) \in G'(r,s,U)(K)$ tels que $U.X$ et $Y.U$ sont de dimensions respectives $ur$ et $us$  (les plus grandes possibles). \\
\noindent On dit que $A$ est bonne si, pour tout entier $d$ compris entre $1$ et $n-1$, les propri\'et\'es suivantes sont satisfaites.\\
\noindent (*) L'alg\`ebre $A$ poss\`ede un bon sous-vectoriel de dimension $d$.\\
\noindent (**) L'action de $\PGL_1(A)$ sur $\G(d,A)$ est g\'en\'eriquement libre.
\end{defi}

A titre d'exercice, le lecteur pourra construire des alg\`ebres qui ne satisfont aucune des deux conditions d\'efinissant une alg\`ebre bonne. L'auteur de ces lignes pense que toute alg\`ebre s\'eparable (i.e. de radical de Jacobson nul) est bonne. Cependant, il n'a pu le montrer que dans le cas d'une alg\`ebre \'etale, ce qui est suffisant pour obtenir de nombreuses applications.

\begin{lem} \label{etalebonne}
Soit $A$ une alg\`ebre \'etale. Alors $A$ est bonne. 
\end{lem}

\begin{dem}
Montrons d'abord que $A$ satisfait la condition (*) de la d\'efinition \ref{bonne}. Soit $n$ la dimension de $A$. Soient $r,s,u$ des entiers positifs tels que $n\geq ru+s$ et $n\geq su+r$. Supposons $s \geq r$; l'autre cas se traite de m\^eme. Soit $t$ un g\'en\'erateur de $A$ comme alg\`ebre (un tel \'el\'ement existe car $K$ est infini). Une base de $A$ comme vectoriel est donc $(1,t,t^2, \ldots, t^{n-1})$. 
 Posons $$X:= <1,t, \ldots, t^{n-r-1}>^\perp \subset A^*,$$ $$Y:= <1,t, \ldots, t^{s-1}> \subset A,$$ $$U:=<1, t^s, t^{2s}, \ldots t^{(u-1) s}> \subset A.$$ On v\'erifie tout de suite que:\\
\noindent i) $Y.U = <1,t, \ldots, t^{su-1}> \subset <1,t, \ldots, t^{n-r-1}>$,\\
\noindent ii)  $(X,Y) \in G'(r,s,U)(K)$.\\
\noindent Pour conclure, il reste donc \`a voir que $U.X$ est de dimension $ur$. Pour ce faire, on peut supposer $r \geq 1$. Par l'absurde, supposons le morphisme canonique $U \otimes X \longrightarrow U.X$ non injectif. On aurait alors une relation de la forme $$\phi_0+ t^s \phi_1 + t^{2s} \phi_2 + \ldots + t^{ms} \phi_m = 0,$$ o\`u les $\phi_i$ appartiennent \`a $X$, avec $0<m<u$ et $\phi_m \neq 0$. Soit $i<m$ un entier positif. Puisqu'on a $n-r-us \geq 0$ et $s \geq r$, $\phi_i$ est nulle sur $<t^{n-r-(m-i)s}, t^{n-r-(m-i)s+1}, \ldots, t^{n-1-(m-i)s}> \subset <1,t, \ldots, t^{n-1-r}>$. Par suite, $t^{is} \phi_i$ est nulle sur $<t^{n-r-ms}, t^{n-r-ms+1}, \ldots, t^{n-1-ms}>$. Donc $t^{ms} \phi_m$ est nulle sur ce sous-espace, ce qui implique que $\phi_m$ est nulle sur 
$<t^{n-r}, t^{n-r+1}, \ldots, t^{n-1}>$. Puisqu'elle appartient \`a $X$, elle est donc identiquement nulle, une contradiction.\\
Montrons maintenant que $A$ satisfait la condition (**) de la d\'efinition \ref{bonne}. Soit $1 \leq r \leq n-1 $ un entier. Montrons que l'action de $\PGL_1(A)$ sur $\G(r,A)$ est g\'en\'eriquement libre. Soit $K':=K(X_0, \ldots, X_{n-1})$ une extension transcendante pure de $K$, et $A':=K'[T] / <T^n+X_{n-1}T^{n-1}+ \ldots+ X_1 T+X_0>$. La propri\'et\'e \`a d\'emontrer \'etant de nature g\'eom\'etrique, et deux alg\`ebres \'etales de m\^eme dimension \'etant g\'eom\'etriquement isomorphes, on peut remplacer $K$ par $K'$ et $A$ par $A'$. On peut donc supposer que $A$ et $K$ sont les seules sous-alg\`ebres de $A$. Il nous faut montrer que le stabilisateur sch\'ematique du point g\'en\'erique de $\G(r,A)$ est trivial. Puisque les points $K$-rationnels de $\G(r,A)$ sont Zariski-denses, il suffit de montrer que le stabilisateur sch\'ematique de tout point $K$-rationnel de $\G(r,A)$ est trivial. Soit donc $E \in \G(r,A)(K)$. Soit $\overline K[\epsilon]$ la $\overline K$-alg\`ebre des nombres duaux. V\'erifions que l'ensemble des  $\overline K[\epsilon]$-points du stabilisateur sch\'ematique de $E$ ne contient que l'\'el\'ement neutre, i.e. que  $$\{ a \in \PGL_1(A)(\overline K[\epsilon]), a.(E\otimes \overline K[\epsilon]) = E\otimes \overline K[\epsilon]\}=\{ 1\}.$$ L'ensemble  $$\{ a \in A, a.E \subset E\}$$ est une sous-alg\`ebre de $A$: c'est donc $K$ ou $A$. On voit facilement que ce dernier cas est exclu. Puisque $\overline K[\epsilon]$ est libre sur $K$, on en d\'eduit  $$ \overline K[\epsilon]=\{ a \in A \otimes \overline K[\epsilon], a.(E\otimes \overline K[\epsilon]) \subset (E\otimes \overline K[\epsilon])\},$$ d'o\`u l'assertion. 
\end{dem}

\begin{lem} \label{bondense}
Soit $A$ une bonne alg\`ebre, de dimension $n$. Les assertions suivantes sont vraies.\\

\noindent i) Soit $u$ un entier compris entre $0$ et $n$. Il existe un ouvert non vide de $\G(u,A)$ dont tout point rationnel est un bon sous-vectoriel.\\
\noindent ii) Soit $0 \leq u_1 < u_2 < \ldots < u_j \leq n$ une collection de $j$ entiers. Soit $\mathbb D_{u_1, \ldots, u_j}(A)$ la vari\'et\'e des drapeaux de $A$, de type $(u_1, \ldots, u_j)$. Il existe un ouvert non vide de $\mathbb D_{u_1, \ldots, u_j}(A)$ dont tout point rationnel $(U_1, \ldots, U_j)$ est form\'e de bons sous-vectoriels.
\end{lem}
\begin{dem}
Montrons i). Soient $r$ et $s$ deux entiers positifs tels que $n\geq ru+s$ et $n\geq su+r$. Il suffit d'exhiber un ouvert non vide $\mathcal U_{r,s} \subset \G(u,A)$ tel que, pour tout $U \in \mathcal U_{r,s}(K)$, il existe $(X,Y) \in G'(r,s,U)(K)$ tel que $U.X$ et $Y.U$ sont de dimensions respectives $ur$ et $us$. En effet, l'intersection des $\mathcal U_{r,s}$, pour $r,s$ comme il pr\'ec\`ede, jouit alors de la propri\'et\'e requise. \\
\noindent Soient donc $r$ et $s$ deux entiers positifs tels que $n\geq ru+s$ et $n\geq su+r$. Soit $U_0$ un bon sous-vectoriel de dimension $u$. Soit $(X,Y) \in G'(r,s,U_0)(K)$ tels que $U_0.X$ et $Y.U_0$ sont de dimensions respectives $ur$ et $us$. Soit $Z \subset A$ un sous-vectoriel, de dimension $ur+s$, tel que $(U_0.X)^\perp \cap Z=Y$. Il existe un ouvert $\mathcal U \subset \G(u,A)$, contenant $U_0$, tel que $U.X$ soit de dimension maximale $ur$ pour tout $U \in \mathcal U( K)$ (de fa\c{c}on explicite, $\mathcal U$ peut s'obtenir comme le lieu o\`u un d\'eterminant ne s'annule pas). Quitte \`a r\'etr\'ecir $\mathcal U$, on peut supposer que $(U.X)^\perp \cap Z$ est de dimension $s=\dim (Z)- \dim(U.X)$ pour tout $U \in \mathcal U(K)$. Pour $U=U_0$, le vectoriel $((U.X)^\perp \cap Z).U$ est \'egal \`a $Y.U_0$, de dimension $su$. R\'etr\'ecissant encore $\mathcal U$, on peut supposer que $((U.X)^\perp \cap Z).U$ est de dimension $su$ pour tout  $U \in \mathcal U(K)$. Pour  $U \in \mathcal U(K)$, posons $Y=(U.X)^\perp \cap Z$. On a alors bien $(X,Y) \in G'(r,s,U)(K)$; de plus $U.X$ et $Y.U$ sont de dimensions respectives $ur$ et $us$. \\
\noindent Montrons ii). Pour $i=1 \ldots j$, soit $\pi_i: \mathbb D_{u_1, \ldots, u_j}(A) \longrightarrow \G(u_i,A)$ le morphisme qui \`a un drapeau $(U_1, \ldots, U_j)$ associe $U_i$. D'apr\`es i), il existe, pour tout $i=1 \ldots j$, un ouvert non vide $\mathcal U_i \subset \G(u_i,A)$ dont tout point rationnel est un bon sous-vectoriel. L'intersection des $\pi_i ^{-1}(\mathcal U_i)$ fait alors l'affaire.
\end{dem}

\begin{defi} \label{defG}
Soit $A$ une bonne alg\`ebre, de dimension $n$. Soit $U \subset A$ un bon sous-vectoriel, de dimension $u$. Soient $r$ et $s$ deux entiers positifs tels que $n\geq ru+s$ et $n\geq su+r$. Les couples $(X,Y) \in G'(r,s,U)$ tels que $U.X$ et $Y.U$ sont de dimensions respectives $ur$ et $us$ forment alors un ouvert, non vide, de $G'(r,s,U)$. Cet ouvert est not\'e $G(r,s,U)$. Sa cl\^oture Zariski, dans $G'(r,s,U)$, est not\'ee $\overline G(r,s,U)$.
\end{defi}

\begin{lem} \label{dom}
Soit $A$ une bonne alg\`ebre, de dimension $n$. Soit $U \subset A$ un bon sous-vectoriel, de dimension $u$. Soient $r$ et $s$ des entiers positifs tels que $n\geq ru+s$ et $n\geq su+r$. Soit $\pi'_1: G'(r,s,U) \longrightarrow \G(r,A^*)$ (resp. $ \pi'_2: G'(r,s,U) \longrightarrow \G(s,A)$) la projection canonique, et $\pi_1$ (resp. $\pi_2$) sa restriction \`a $G(r,s,U)$. Les propri\'et\'es suivantes sont vraies.\\

\noindent i) $G(r,s,U)$ est munie d'une action naturelle de $\PGL_1(A)$. \\
\noindent ii) Les morphismes $\pi_1$ et $\pi_2$ sont dominants. \\
\noindent iii)  $G(r,s,U)$ est une vari\'et\'e g\'eom\'etriquement connexe et lisse.\\
\noindent Sa dimension est $r(n-r)+s(n-s)-sru$.\\
\noindent iv) Le ferm\'e $\overline G(r,s,U) \subset G'(r,s,U)$ est l'unique composante irr\'eductible $C' \subset G'(r,s,U)$ telle que la restriction de $\pi'_1$ et de $\pi'_2$ \`a $C'$ est dominante.
\end{lem}
\begin{dem}
Pour i), on v\'erifie sans difficul\'e que l'ouvert $G(r,s,U) \subset G'(r,s,U)$ est stable par l'action de $\PGL_1(A)$.\\
\noindent Pour d\'emontrer les autres points, on peut supposer $K$ alg\'ebriquement clos.\\
\noindent La d\'efinition de $G'(r,s,U)$ montre que cette vari\'et\'e est localement d\'efinie par $rsu$ \'equations. Par suite, toutes ses composantes irr\'eductibles (donc aussi celles de $G(r,s,U)$) sont de dimension $\geq s(n-s)+r(n-r)-rsu$. Nous allons montrer dans un instant que l'espace tangent en tout point $K$-rationnel de $G(r,s,U)$ est de dimension $r(n-r)+s(n-s)-rsu$. Par ce qui vient d'\^etre dit, cela implique que $G(r,s,U)$ est lisse et que toutes ses composantes irr\'eductibles sont de dimension $r(n-r)+s(n-s)-rsu$. Soit $C$ une telle composante. Par le lemme 2.4 ii), la fibre de ${\pi_1}_{|C}$ en $(X,Y) \in C(K)$ est un ouvert non vide de $\G(s,(UX)^\perp)$, donc de dimension $s(n-ru-s)$. Un calcul imm\'ediat montre que $s(n-ru-s)+\dim \G(r,A^*)=\dim C$; ${\pi_1}_{|C}$ est donc dominante. Vu que la fibre de $\pi_1$ en tout point $K$-rationnel est irr\'eductible (lemme 2.4 ii)), ceci n'est possible que si $C$ est la seule composante irr\'eductible de $G(r,s,U)$. Ceci d\'emontre ii) et iii).  L'assertion d'unicit\'e du point iv) r\'esulte alors de ce que les fibres g\'eom\'etriques de $\pi'_1$ (resp. de $\pi'_2$) sont irr\'eductibles. \\
Soit donc $(X,Y) \in G(r,s,U) ( K)$.  On sait que l'espace tangent en $X$ (resp. $Y$) de $\G(r,A^*)$ (resp. de $\G(s,A)$) s'identifie canoniquement \`a $\Hom(X,A^* /X)$ (resp. \`a $\Hom(Y,A /Y)$). En diff\'erenciant les \'equations d\'efinissant $G'(r,s,U)$, on voit que son espace tangent en $(X,Y) \in G(r,s,U)(K)$ n'est autre que le noyau du morphisme $$\Theta: \Hom(X,A^* /X) \oplus \Hom(Y,A/Y) \longrightarrow (X \otimes Y \otimes U)^*,$$ $$(f,g) \mapsto (x \otimes y \otimes z  \mapsto f(x)(yz) + x(g(y)z)). $$ Ce morphisme est surjectif. En effet, soit $\phi \in X^*$ et $\psi \in (Y \otimes U)^*$. Puisque le morphisme `multiplication' $  Y \otimes U \longrightarrow Y.U$ est un isomorphisme, il existe $\psi' \in A^*$ tel que $\psi'(yz)= \psi(y \otimes z)$ pour tous $y \in Y$, $z \in U$. Soit $f: X \longrightarrow A^*/X$ le morphisme envoyant $x \in X$ sur la classe de $\phi(x) \psi'$. On a $\Theta(f,0)=\phi \otimes \psi$. Puisque $\phi$ et $\psi$ sont arbitraires, cela montre la surjectivit\'e de $\Theta$. L'espace tangent en $(X,Y)$ \`a  $G(r,s,U)$ est donc de dimension $r(n-r)+s(n-s)-sru$, cqfd.\end{dem}

Pour illustrer l'utilisation que nous ferons des bons sous-vectoriels, commen\c{c}ons par un petit lemme.
\begin{lem} \label{AetAstar}
 Soit $A$ une bonne alg\`ebre, de dimension $n$. Soit $n=qd$ une d\'ecomposition de $n$ en produit de deux entiers. Soit $U \subset A$ un bon sous-vectoriel, de dimension $q-1$. Alors l'application rationnelle $\PGL_1(A)$-\'equivariante $$\Phi_U: \G(d,A) \longrightarrow \G(d,A^*), $$   $$ Y \mapsto (Y.U)^\perp,$$ est birationnelle.
\end{lem}
\begin{dem}
Soient $\pi_1: G(d,d,U) \longrightarrow \G(d,A^*)$ et $\pi_2: G(d,d,U) \longrightarrow \G(d,A)$ les projections canoniques, dominantes par le lemme \ref{dom} ii). D'apr\`es le lemme \ref{lemG'} ii), la fibre g\'en\'erique de $\pi_1$ est un point (grassmannienne des sous-espaces de dimension $d$ d'un vectoriel de dimension $n-(q-1)d=d$), donc $\pi_1$ est birationnelle puisque sa source est int\`egre. On peut exhiber l'inverse birationnel de $\pi_1$ comme $$X \in \G(d,A^*) \mapsto (X,(U.X)^\perp) \in G(d,d,U).$$ De m\^eme, $\pi_2$ est birationnelle. Puisque $\Phi_U =\pi_1 \circ \pi_2 ^{-1}$, on a l'assertion voulue.
\end{dem}

Nous allons maintenant d\'evelopper quelques r\'esultats plus raffin\'es.

\begin{lem} \label{biendef}
Les objets et notations sont ceux de la d\'efinition \ref{defG}. Distinguons deux cas.\\
\noindent i)  Supposons $r \geq s$. Ecrivons la division euclidienne $r=sq+t$. Soit $U' \subset A$ un bon sous-vectoriel, contenant $U$, de dimension $u+q$. Il existe alors un ouvert non vide $\mathcal U \subset G(r,s,U)$ tel que, pour tout $(X,Y) \in  \mathcal U (\overline K)$, les deux conditions suivantes sont v\'erifi\'ees.\\
\noindent a) $\dim(Y.U')=s(u+q)$,\\
\noindent b) $\dim (X \cap (Y.U')^\perp)=t$.\\

\noindent ii) Supposons $s \geq r$. Ecrivons la division euclidienne $s=rq+t$. Soit $U' \subset A$ un bon sous-vectoriel, contenant $U$, de dimension $u+q$. Il existe alors un ouvert non vide $\mathcal U \subset G(r,s,U)$ tel que, pour tout $(X,Y) \in  \mathcal U (\overline K)$, les deux conditions suivantes sont v\'erifi\'ees.\\
\noindent a) $\dim(U'.X)=r(u+q)$,\\
\noindent b) $\dim (Y \cap (U'.X)^\perp)=t$.\\
\end{lem}
 \begin{dem}
Donnons la d\'emonstration de i); celle de ii) est identique. Les conditions a) et b) d\'efinissent un sous-sch\'ema ouvert $\mathcal U \subset G(r,s,U)$, pour la raison suivante. Supposons satisfaite la condition a), qui est ouverte (elle exprime que le morphisme lin\'eaire `multiplication' $Y \otimes U' \longrightarrow A$ est injectif). Soient $A$ et $B$ les vectoriels $X$ et $(Y.U')^\perp$, vus comme sous-vectoriels de $C:=(Y.U)^\perp$. La condition b) exprime alors simplement que $A ^\perp \subset C^*$ et $B ^\perp \subset C^*$ (de dimensions respectives $n-su-r$ et $n-su-(n-s(u+q))=sq$) sont en somme directe. C'est aussi une condition ouverte. Il nous faut donc montrer que $\mathcal U$ est non vide. Soit $ \pi: G(r,s,U) \longrightarrow \G(s,A)$ la projection, dominante par le lemme \ref{dom}. Il existe donc $(X',Y) \in G(r,s,U)(K)$ satisfaisant a). La fibre de $\pi$ en $(X',Y)$ est un ouvert non vide de $\G(r,(Y.U)^\perp)$. L'ouvert form\'e des $X \in \G(r,(Y.U)^\perp)$ satisfaisant  b) \'etant aussi non vide, n'importe quel point rationnel $X$ de l'intersection de ces deux ouverts est tel que $(X,Y)$ satisfait a) et b), cqfd.
\end{dem}

\begin{lem} \label{valeurs}
Soit $A$ une bonne alg\`ebre, de dimension $n$. Soient $r,s,u$ trois entiers positifs, tels que $n \geq su +r$ et $n \geq ru +s$.\\
\noindent  Supposons d'abord $r \geq s$. Ecrivons la division euclidienne $r=qs+t$. Soient $U \subset U' \subset A$ deux bons sous-vectoriels, de dimensions $u$ et $u+q$. On a $$ t(u+q)+s \leq s(u+q)+t = su+r \leq n;$$ on peut donc consid\'erer la vari\'et\'e $G'(t,s,U')$. D'apr\`es le lemme \ref{biendef}, on a une  application rationnelle bien d\'efinie $$\phi: G(r,s,U) \longrightarrow G'(t,s,U'),$$  au moyen de la formule $$ (X,Y) \mapsto (X \cap (Y.U')^\perp, Y).$$ Cette application rationnelle est \`a valeurs dans $\overline G(t,s,U')$.\\

 Si $s > r$,  \'ecrivons la division euclidienne $s=rq+t$. Soient $U \subset U' \subset A$ deux bons sous-vectoriels, de dimensions $u$ et $u+q$.  On d\'efinit alors de mani\`ere semblable une application rationnelle $$\phi: G(r,s,U) \longrightarrow G'(r,t,U'),$$  au moyen de la formule $$ (X,Y) \mapsto (X, Y \cap (U'.X)^\perp).$$
Cette application rationnelle est \`a valeurs dans $\overline G(r,t,U')$.
 \end{lem}

\begin{dem}
Supposons $r \geq s$; l'autre cas se traite de m\^eme. Puisque  $G(r,s,U)$ est irr\'eductible (lemme \ref{dom} iii)), $\phi$ prend ses valeurs dans une composante irr\'eductible $C \subset  G'(t,s,U')$. Soient $\pi:G(r,s,U) \longrightarrow \G(s,A)$ et $\pi': G'(t,s,U') \longrightarrow \G(s,A)$ les projections. On a $\pi' \circ \phi = \pi$. Puisque $\pi$ est dominante (lemme \ref{dom} ii)), on en d\'eduit que la restriction de $\pi'$ \`a $C$ est dominante. D'apr\`es le lemme \ref{dom} iv), on en tire la conclusion souhait\'ee.
\end{dem}

\begin{prop} \label{propphi}
Soit $A$ une bonne alg\`ebre, de dimension $n$. Soient $r,s,u$ trois entiers positifs, tels que $n \geq su +r$ et $n \geq ru +s$.\\
\noindent  Supposons d'abord $r \geq s$. Ecrivons la division euclidienne $r=qs+t$.  Soient $U \subset U' \subset A$ deux bons sous-vectoriels, de dimensions $u$ et $u+q$. Soit $$\phi: G(r,s,U) \longrightarrow \overline G(t,s,U')$$ l'application rationnelle du lemme \ref{valeurs}. Les propositions suivantes sont vraies.\\

\noindent i) $\phi$ est $\PGL_1(A)$- \'equivariante. \\
\noindent ii) Soit $\mathcal U$ l'ouvert du lemme \ref{biendef}. Soit $L$ une extension de $K$. La fibre du morphisme $\phi_{|\mathcal U} $ en $(X, Y) \in G(t,s,U')(L)$ est un ouvert de $\G_L(sq,  (Y.U)^{\perp}/X)$.\\
\noindent iii) $\phi$ est dominante.\\

Si $s > r$, \'ecrivons la division euclidienne $s=qr+t$. Soient $U \subset U' \subset A$ deux bons sous-vectoriels, de dimensions $u$ et $u+q$.  Soit $$\phi: G(r,s,U) \longrightarrow \overline G(r,t,U')$$ l'application rationnelle du lemme \ref{valeurs}. On a alors les assertions analogues \`a i), ii) et iii).

\end{prop}

\begin{dem}
La premi\`ere assertion est une v\'erification facile. Montrons ii). La fibre de $\phi_{|\mathcal U} $ en $(X,Y)$ est un ouvert du sous-sch\'ema ferm\'e des $X' \in \G_{ L}(r,A^*)$  v\'erifiant $X \subset X' \subset (Y.U)^{\perp}$, qui n'est autre que $\G_{ L}(sq,  (Y.U)^{\perp}/X)$. Cette fibre est donc soit vide, soit de dimension  $sq(n-us-t-sq)=sq(n-us-r)$. D'apr\`es le calcul de dimension fait au point ii) du lemme \ref{dom}, on en d\'eduit que la cl\^oture de l'image de $\phi$ est de dimension $$\dim\mathcal U-sq(n-us-r)=s(n-s)+r(n-su-r)-sq(n-us-r)$$ $$ =s(n-s)+t(n-su-r)=s(n-s)+t(n-s(u+q)-t)=\dim \overline G(t,s,U').$$ Ceci d\'emontre iii). 
\end{dem}

\begin{defi} \label{defiphi}
Soit $A$ une bonne alg\`ebre, de dimension $n$. Soient $r,s,u$ trois entiers positifs, tels que $n \geq su +r$ et $n \geq ru +s$.\\
\noindent  Supposons d'abord $r \geq s$. Ecrivons la division euclidienne $r=qs+t$.  Soient $U \subset U' \subset A$ deux bons sous-vectoriels, de dimensions $u$ et $u+q$. Le lemme \ref{propphi} permet de d\'efinir une application rationnelle, dominante, $\PGL_1(A)$- \'equivariante, \`a fibres rationnelles $$\phi_{U,U'}: G(r,s,U) \longrightarrow  G(t,s,U').$$

\noindent  Supposons maintenant $s > r$. Ecrivons la division euclidienne $s=qr+t$.  Soient $U \subset U' \subset A$ deux bons sous-vectoriels, de dimensions $u$ et $u+q$. On d\'efinit de la m\^eme mani\`ere une application rationnelle, dominante, $\PGL_1(A)$- \'equivariante, \`a fibres rationnelles $$\phi_{U,U'}: G(r,s,U) \longrightarrow  G(r,t,U').$$
\end{defi}

\begin{rem} \label{iterer}
On peut it\'erer la construction pr\'ec\'edente, en suivant l'algorithme d'Euclide, et composer les applications rationnelles $\phi_{U,U'}$  entre elles. De fa\c{c}on pr\'ecise, cela est l'objet de la prochaine section, o\`u nous d\'emontrons le th\'eor\`eme principal.
\end{rem}

\section{Le th\'eor\`eme principal}

\begin{defi}\label{defiprinc}
Soit $A$ une bonne alg\`ebre, de dimension $n$. Soit $1<r<n$ un entier. Soit $$n=r_0, r=r_1,r_2, \ldots, r_s=\mathrm{pgcd}(n,r),r_{s+1}=0$$ la suite d'entiers d\'efinie selon l'algorithme d'Euclide, en prenant pour $r_{i+1}$ le reste de la division euclidienne de $r_{i-1}$ par $r_i$. Soient $q_i$ les entiers tels que $r_{i-1}=q_i r_i + r_{i+1}$. Soit $$\mathcal U = (U_0=\{0\} \subset U_1 \subset U_2 \subset ... \subset U_s)$$ un drapeau form\'e de bons sous-vectoriels de $A$, tels que $\ddim(U_i)-\ddim(U_{i-1})=q_i$ (un tel drapeau existe par le lemme \ref{bondense}). \\
\noindent Supposons $s$ impair. D'apr\`es la d\'efinition \ref{defiphi}, il existe alors une suite d'applications rationnelles dominantes, $\PGL_1(A)$-\'equivariantes, \`a fibres rationnelles: $$ \phi_{U_0,U_1}: G(r_0,r_1,U_0) \longrightarrow G(r_2,r_1,U_1),$$  $$ \phi_{U_1,U_2}: G(r_2,r_1,U_1) \longrightarrow G(r_2,r_3,U_2),$$   $$ \phi_{U_2,U_3}: G(r_2,r_3,U_2) \longrightarrow G(r_4,r_3,U_3),$$ $$ \ldots $$  $$ \phi_{U_{s-1},U_s}: G(r_{s-1},r_s,U_{s-1}) \longrightarrow G(r_{s+1},r_s,U_s).$$ Puisque $G(r_0,r_1,U_0)$ (resp. $G(r_{s+1},r_s,U_s)$) est un ouvert de $\G(r,A)$ (resp. $\G(\mathrm{pgcd}(r,n),A)$), on peut voir la compos\'ee de toutes ces applications comme une application rationnelle, dominante, $\PGL_1(A)$-\'equivariante, \`a fibres rationnelles $$ \Phi_ {\mathcal U}: \G(r,A) \longrightarrow \G(\mathrm{pgcd}(r,n),A).$$ Si $s$ est pair, le m\^eme proc\'ed\'e fournit une application rationnelle, dominante, $\PGL_1(A)$-\'equivariante, \`a fibres rationnelles $$ \Phi_{\mathcal U}: \G(r,A) \longrightarrow \G(\mathrm{pgcd}(r,n),A^*).$$ 

\end{defi}

\begin{prop}\label{propetale}
Soit $A$ une  bonne alg\`ebre, de dimension $n$. Soient $r,s,u$ trois entiers positifs et non nuls, tels que  $n \geq ru +s$ et  $n \geq su +r$.\\
\noindent Supposons $r \geq s$. Ecrivons la division euclidienne $r=qs+t$. Soit $U \subset U' \subset A$ deux bons sous-vectoriels, de dimensions $u$ et $u+q$. Soit $$\phi_{U,U'}: G(r,s,U) \longrightarrow G(t,s,U')$$ l'application rationnelle, dominante et $\PGL_1(A)$-\'equivariante de la d\'efinition \ref{defiphi}. Les propositions suivantes sont vraies.\\

\noindent i) L'action de $\PGL_1(A)$ sur la source et le but de $\phi_{U,U'}$ est g\'en\'eriquement libre.\\

\noindent Soit $ \tilde G(r,s,U)$ (resp. $\tilde G(t,s,U')$) un ouvert de $G(r,s,U)$ (resp. $G(t,s,U')$), bon pour l'action de $\PGL_1(A)$, et tel que $\phi_{U,U'}$ induise un morphisme partout d\'efini $$  \phi: \tilde G(r,s,U) \longrightarrow \tilde G(t,s,U').$$

\noindent ii) Le morphisme quotient $$ \overline \phi: \tilde G(r,s,U)/\PGL_1(A) \longrightarrow \tilde G(t,s,U')/ \PGL_1(A)$$ induit une extension transcendante pure sur les corps de fonctions.\\

\noindent Si  $s > r$, \'ecrivons la division euclidienne $s=qr+t$. Soit $U \subset U' \subset A$ deux bons sous-vectoriels, de dimensions $u$ et $u+q$.  Soit $$\phi_{U,U'}: G(r,s,U) \longrightarrow G(r,t,U')$$ l'application rationnelle, dominante et $\PGL_1(A)$-\'equivariante de la d\'efinition \ref{defiphi}. On a alors les assertions analogues \`a i) et ii).

 \end{prop}
\begin{dem}
On suppose $r \geq s$; l'autre cas se traite de m\^eme.\\
\noindent i) Soit $d=\mathrm{pgcd}(r,s)$. On peut it\'erer la construction de la d\'efinition \ref{defiphi} pour construire une application rationnelle, dominante, $\PGL_1(A)$-\'equivariante $$\theta: G(r,s,U) \longrightarrow \G(d,A)$$ ou bien  $$\theta: G(r,s,U) \longrightarrow \G(d,A^*)=\G(n-d,A).$$ Par un argument \'el\'ementaire laiss\'e au lecteur, ceci montre que l'action de  $\PGL_1(A)$ sur $G(r,s,U)$ est g\'en\'eriquement libre si elle l'est sur $\G(d,A)$ (ou sur $\G(n-d,A)$), ce qui est le cas par d\'efinition d'une alg\`ebre  bonne.\\

\noindent Montrons ii). Soit $\pi:\tilde G(t,s,U') \longrightarrow \tilde G(t,s,U')/ \PGL_1(A)$ la projection canonique. Soit $L$ une extension de $K$, et $y \in (\tilde G(t,s,U')/ \PGL_1(A))(L)$. Montrons que la fibre de $\overline \phi$ en $y$ est $L$-rationnelle. Soit $ \overline L$ une cl\^oture s\'eparable de $L$. Soit $x=(X,Y) \in \tilde G(t,s,U')(\overline L)$ s'envoyant sur $y$ par $\pi$ (un tel point existe car $\PGL_1(A)$ est lisse sur $K$). De la fa\c{c}on habituelle, on peut associer \`a $x$ un 1-cocycle galoisien $$c_\sigma: Gal(\overline L /L) \longrightarrow \PGL_1(A)(\overline L).$$ Il d\'efinit donc une alg\`ebre simple centrale $\mathcal A$, isomorphe au tordu de $\End_L(A \otimes L)$ par le $\PGL_1(A)$-torseur `fibre de $\pi$ en $y$'. Si $E$ est un sous-$\overline L$-vectoriel de $A \otimes \overline L$ ou de $A^* \otimes \overline L$, et $\sigma \in  Gal(\overline L /L)$, on pose  $$\sigma * E=c_\sigma (\sigma(E)).$$ L'op\'eration $*$  laisse $\P_{\overline L}(X) \subset \P_{\overline L}(A^* \otimes \overline L)$ et $\P_{\overline L }(Y) \subset \P_{\overline L}(A \otimes \overline L)$ stables. Par suite, les vari\'et\'es $SB(t,\mathcal A^{op})$ et $SB(s,\mathcal A)$ ont des points $L$-rationnels (correspondant \`a $X$ et $Y$ respectivement). Par la proposition \ref{propphi}, la fibre de $\phi$ en $x$ est un ouvert de $\G_{\overline L}(sq,(Y.U)^\perp /X)$. Puisque $(Y.U)^\perp$ et $X$ sont fixes par $*$, l'op\'eration $*$ induit une action semi-lin\'eaire de $Gal(\overline L /L)$ sur $\G_{\overline L}(sq,(Y.U)^\perp /X)$, encore not\'ee  $*$. Par ce qui pr\'ec\`ede, la fibre de $\overline \phi$ en $y$ est un ouvert de la $L$-forme de $\G_{\overline L}(sq,(Y.U)^\perp /X)$ donn\'ee par l'action $*$, qui est une $L$-vari\'et\'e de Severi-Brauer g\'en\'eralis\'ee. Une telle vari\'et\'e \'etant rationnelle d\`es qu'elle poss\`ede un point rationnel, il nous faut voir que $(Y.U)^\perp /X$ poss\`ede des sous-$\overline L$-espaces vectoriels stables par $*$, de dimension $sq$. Consid\'erons la condition suivante, portant sur un sous-espace $Z \in \G(s(u+q),A^*)(\overline L)$: \\
\noindent (C): $(Y.U)^\perp \cap Z$  est de dimension $sq$, et intersecte trivialement $X$.\\
\noindent Elle d\'efinit un ouvert de $\G_{\overline L}(s(u+q),A^* \otimes \overline L)$. Cet ouvert est non vide. En effet, on a $s(q+u)+t \leq n$, donc $sq+t \leq n-su$; autrement dit $sq+\dim X \leq \dim (Y.U)^\perp$. Il existe donc un sous-espace $Z' \subset (Y.U)^\perp$, de dimension $sq$, et d'intersection triviale avec $X$. Un tel $Z'$ s'\'ecrit toujours comme $(Y.U)^\perp \cap Z$ pour un  $Z \in \G(s(u+q),A^*)(\overline L)$, puisque la codimension de $(Y.U)^\perp$ dans $A^*$ est $su$. D'o\`u la non-vacuit\'e annonc\'ee. Maintenant, l'entier $s(u+q)$ est un multiple de $s$; les points $L$-rationnels de $SB(s(u+q),\mathcal A)$ sont donc Zariski-denses. Par suite, il existe un sous-$\overline L$-vectoriel $Z   \in \G(s(u+q),A^*)(\overline L)$, $*$-stable, satisfaisant (C). L'image de $(Y.U)^\perp \cap Z$ par la projection $(Y.U)^\perp \longrightarrow (Y.U)^\perp /X$ fait alors l'affaire.\\
\noindent On applique ceci lorsque $L$ est le corps des fonctions de $\tilde G(t,s,U')/ \PGL_1(A)$ et $y$ son point g\'en\'erique, ce qui d\'emontre ii).

\end{dem}

Nous disposons maintenant du mat\'eriel requis pour d\'emontrer le th\'eor\`eme principal de cette section.

\begin{thm} \label{thprinc}
 Soit $A$ une  bonne alg\`ebre, de dimension $n$.  Soit $0<r<n$ un entier. Soit $d=\mathrm{pgcd}(r,n)$. Alors la vari\'et\'e $\G(r,A)$ est birationnelle, de mani\`ere $\PGL_1(A)$-\'equivariante, au produit de $\G(d,A)$ par un espace affine sur lequel $\PGL_1(A)$ agit trivialement. 
\end{thm}
\begin{dem}
En combinant la d\'efinition \ref{defiprinc} et la proposition \ref{AetAstar}, on dispose d'une application rationnelle, dominante, $\PGL_1(A)$- \'equivariante $$ \Phi: \G(r,A) \longrightarrow \G(d,A)$$ Soit $\tilde{\G}(r,A)$ (resp. $\tilde{\G}(d,A)$) un bon ouvert de $\G(r,A)$ (resp. $\G(d,A)$) pour l'action de $\PGL_1(A)$, tels que  $\Phi$ induise un morphisme partout d\'efini $$ \Phi': \tilde{ \G}(r,A) \longrightarrow \tilde{\G}(d,A).$$  Par applications it\'er\'ees du point iv) de la proposition \ref{propetale}, on voit que le morphisme quotient  $$ \overline \Phi': \tilde{\G}(r,A)/\PGL_1(A) \longrightarrow \tilde{\G}(d,A)/ \PGL_1(A) $$ induit une extension transcendante pure sur les corps de fonctions. Or le morphisme canonique $$  \tilde{ \G}(r,A)  \longrightarrow  \tilde{ \G}(d,A)  \times_{\tilde{\G}(d,A) /\PGL_1(A)  }( \tilde{\G}(r,A)/\PGL_1(A))$$ induit par le diagramme $$\xymatrix{\tilde{ \G}(r,A) \ar[d]\ar[r] & \tilde{ \G}(d,A) \ar[d]\\ \tilde{\G}(r,A)/\PGL_1(A) \ar[r] &  \tilde{\G}(d,A)/\PGL_1(A)}$$est un isomorphisme (un morphisme entre deux torseurs est un isomorphisme!), ce qui termine la d\'emonstration.
\end{dem}

\begin{coro}
Sous les hypoth\`eses du th\'eor\`eme \ref{thprinc}, supposons de plus que $d=1$. Alors la vari\'et\'e  $\G(r,A) / \PGL_1(A)$ (d\'efinie \`a isomorphisme birationnel pr\`es) est rationnelle.
\end{coro}
\begin{dem}
Appliquant le th\'eor\`eme \ref{thprinc}, on voit que  $\G(r,A)$ est birationnelle, de mani\`ere $\PGL_1(A)$-\'equivariante, au produit de $\P(A)$ par un espace affine. Il revient au m\^eme de dire qu'elle est birationnelle, de mani\`ere $\PGL_1(A)$-\'equivariante, au produit de $\PGL_1(A)$ par un espace affine, d'o\`u l'assertion.
\end{dem}
\begin{rem}
Soit $V$ un vectoriel de dimension $n$, et $A= \End(V)$. Le lecteur pourra montrer, \`a titre d'exercice, que le quotient birationnel $\G(n,A)/ \PGL_1(A)$ est stablement birationnel \`a l'espace classifiant de $\PGL_1(A)=\PGL(V)$. Sa (stable) rationalit\'e est un probl\`eme majeur et tr\`es largement ouvert. Elle n'est en effet connue que lorsque $n$ divise $420$ (\cite{BLB}).
\end{rem}

\section{Quelques applications du th\'eor\`eme \ref{thprinc}}

Nous \'enon\c{c}ons dans cette section, sous forme de propositions, plusieurs corollaires du th\'eor\`eme principal, qui s'en d\'eduisent par torsion. Pour des rappels sur cet important proc\'ed\'e, nous renvoyons le lecteur \`a \cite{F}, proposition 2.12 et lemme 2.14.

\begin{prop}\label{SB}
Soit $0<r<n$ deux entiers, et $d=\mathrm{pgcd}(r,n)$. Soit $A$ une alg\`ebre simple centrale de degr\'e $n$. Alors la vari\'et\'e de Severi-Brauer g\'en\'eralis\'ee $SB(r,A)$ est birationnelle au produit de  $SB(d,A)$ par un espace affine de dimension convenable.
\end{prop}
\begin{dem}
Soit $L$ une sous-alg\`ebre \'etale commutative maximale de $A$. Il est classique (suite exacte longue de cohomologie combin\'ee au th\'eor\`eme 90 de Hilbert) que la cohomologie du groupe alg\'ebrique $T:=R_{L/K}(\Gm)/ \Gm=\PGL_1(L)$ s'identifie au noyau de la fl\`eche naturelle $Br(K) \longrightarrow Br(L)$; il correspond donc \`a $A$ un $T$-torseur $P$ tel que $SB(m,A)$ s'identifie au tordu de $\G(m,L)$ par $P$, pour tout $m=1 \ldots n$. On d\'eduit tout de suite le corollaire par torsion \`a partir du th\'eor\`eme \ref{thprinc}, appliqu\'e \`a $L$.
\end{dem}

\begin{defi}

Soit $G$ un groupe alg\'ebrique agissant sur une vari\'et\'e g\'eom\'etriquement int\`egre $X$. Nous dirons que l'action de $G$ sur $X$ est $\textit{parfaite}$ si $X$ est birationnelle, de mani\`ere $G$-\'equivariante, au produit de $G$ par un espace affine sur lequel $G$ agit trivialement.
\end{defi}
\begin{rem}
Soit $G$ un groupe alg\'ebrique agissant sur une vari\'et\'e g\'eom\'etriquement int\`egre $X$. Le lecteur se convaincra ais\'ement que l'action de $G$ sur $X$ est parfaite si et seulement si les conditions suivantes sont v\'erifi\'ees.\\
\noindent i) L'action de $G$ sur $X$ est g\'en\'eriquement libre.\\
\noindent ii) Soit $U \subset X$ un ouvert bon pour l'action de $G$. Alors le $G$-torseur $U \longrightarrow U/G$ poss\`ede des sections rationnelles, et sa base est rationnelle.
\end{rem}

\begin{prop} \label{actionparfaite}
Soit $V$ un vectoriel de dimension $r$, et $A$ une bonne alg\`ebre de dimension $n$. Supposons $r \leq n$ et $pgcd(n,r)=1$. Soit $G$ le conoyau du morphisme $$ \G_m \longrightarrow \GL(V) \times \GL_1(A),$$ $$ x \mapsto (x,x^{-1}).$$ Le groupe $G$ agit sur $\A(V \otimes A)$ par la formule $$(f,x).(v \otimes a)=f(v) \otimes xa,$$ pour $f \in \GL(V)$, $x \in \GL_1(A)$, $v \in V$ et $a \in A$.\\
\noindent Cette action est parfaite.
\end{prop}
\begin{dem}
Soit $U$ l'ouvert de $\A(V \otimes A)= \A( \Hom(V^*,A))$ form\'e des applications lin\'eaires de rang maximal $r$. L'ouvert $U$ est stable par l'action naturelle de $\GL(V)$ sur $\A(V \otimes A)$, et le quotient $U/ \GL(V)$ n'est autre que $\G(r,A)$. La projection naturelle $U \longrightarrow \G(r,A)$ est un $\GL(V)$-torseur. Puisque l'action de $\PGL_1(A)=G/\GL(V)$ sur  $\G(r,A)$ est g\'en\'eriquement libre,  on en d\'eduit que l'action de $G$ sur $\A(V \otimes A)$ l'est aussi. On sait alors (d\'efinition \ref{bonouvert}) qu'il existe un ouvert $U' \subset \A(V \otimes A)$, bon pour l'action de $G$. Par le th\'eor\`eme 90 de Hilbert, la projection $U' \longrightarrow U'/\GL(V)$ poss\`ede des sections rationnelles. D'apr\`es le th\'eor\`eme \ref{thprinc}, on sait que l'action de  $\PGL_1(A)$ sur $\G(r,A)$ est parfaite. Ceci se traduit ici en disant que le $\PGL_1(A)$-torseur $U'/ \GL(V) \longrightarrow U'/G$ poss\`ede des sections rationnelles, et que sa base est rationnelle. En d\'efinitive, on voit que le $G$-torseur $U' \longrightarrow U'/G$ poss\`ede des sections rationnelles, et que sa base est rationnelle. En d'autres termes, l'action de $G$ sur $\A(V \otimes A)$ est parfaite.
\end{dem}

\begin{prop}\label{tens}
Soient $A$ et $B$ deux alg\`ebres simples centrales, de degr\'es premiers entre eux. Alors $SB(A\otimes B)$ est birationnelle au produit de $SB(A) \times SB(B)$ par un espace affine de dimension convenable. 
\end{prop}
\begin{dem}
Soit $n$ (resp. $r$) le degr\'e de $A$ (resp. $B$). On suppose $r<n$.
Soit $L$ (resp. $M$) une sous-alg\`ebre \'etale commutative maximale de $A$ (resp. $B$).  Soit $G$ le conoyau du morphisme $$ \G_m \longrightarrow \GL(M) \times \GL_1(L),$$ $$ x \mapsto (x,x^{-1}).$$ Soit $G'$ le conoyau du morphisme $$ \G_m \longrightarrow \GL_1(M) \times \GL_1(L),$$ $$ x \mapsto (x,x^{-1});$$ c'est un sous-groupe de $G$. Par le th\'eor\`eme 90 de Hilbert, le torseur $\GL(M) \longrightarrow  \GL_1(M) \backslash \GL(M)$ poss\`ede des sections rationnelles. Sa base \'etant rationnelle (exercice laiss\'e au lecteur), il suit que  l'action de $\GL_1(M)$ sur $\GL(M)$ (par multiplication \`a gauche) est parfaite, donc que l'action de $G'$ sur $G$ est parfaite. Puisque l'action de $G$ sur $\A(M \otimes L)$ est parfaite par la proposition \ref{actionparfaite}, on en d\'eduit -et c'est nettement plus faible- que l'action de $G'$ sur $\A(L \otimes M)$ est parfaite. Passant au quotient par $\G_m$, on voit finalement que l'action de  $\PGL_1(L) \times \PGL_1(M)$ sur $\P(L \otimes M)$ est parfaite. En d'autres termes, la vari\'et\'e $\P(L \otimes M)$ est birationnelle, de mani\`ere $\PGL_1(L) \times \PGL_1(M)$-\'equivariante, au produit de $\P(L) \times \P(M)$ par un espace affine. Le r\'esultat annonc\'e suit par torsion, comme dans la preuve de la proposition \ref{SB}.
\end{dem}

\begin{prop}(\cite{Kly}) \label{tores}
Soient $L$ et $M$ deux alg\`ebres \'etales, de degr\'es premiers entre eux. Soit $T$ le conoyau du morphisme $$ \GL_1
(M) \times \GL_1(L) \longrightarrow \GL_1 ( M \otimes L),$$ $$(x,y) \mapsto x \otimes y;$$ c'est un tore alg\'ebrique. La vari\'et\'e $T$ est rationnelle.
\end{prop}
\begin{dem}
Soit $G'$ le conoyau du morphisme $$ \G_m \longrightarrow \GL_1(M) \times \GL_1(L),$$ $$ x \mapsto (x,x^{-1}).$$ Dans la preuve de la proposition \ref{tens}, on a vu que l'action de $G'$ sur $\A(M \otimes L)$ -donc aussi sur $\GL_1 ( M \otimes L)$- est parfaite. Le r\'esultat suit par d\'efinition d'une action parfaite, vu que $T$ n'est autre que $\GL_1 ( M \otimes L)/G'$.
\end{dem}

\begin{prop}(\cite{Kra}, Theorem 3.1)
Soit $A$ une alg\`ebre simple centrale de degr\'e $n$. Soit $0 <r <n$ un entier. Alors les vari\'et\'es $SB(r,A)$ et $SB(r,A^{op})$ sont birationnellement isomorphes.
\end{prop}
\begin{dem}
Au vu de la proposition \ref{SB}, on peut supposer que $r$ divise $n$. Soient $L$ et $P$ comme dans la preuve de cette proposition.  D'apr\`es le lemme \ref{AetAstar},  il existe un isomorphisme birationnel et $\PGL_1(L)$-\'equivariant $$f: \G(r,L)  \longrightarrow \G(r,L^*).$$  En tordant $f$ par $P$, on obtient un isomorphisme birationnel entre $SB(r,A)$ et $SB(r,A^{op})$, cqfd.
\end{dem}

\begin{prop}
Soit $r,n$ deux entiers premiers entre eux, v\'erifiant $3 \leq r \leq n-3$. Soient $A, A'$ deux alg\`ebres simples centrales de degr\'e $n$, engendrant le m\^eme sous-groupe du groupe de Brauer de $K$.  Alors les vari\'et\'es $SB(r,A)$ et $SB(r,A')$ sont birationnelles.
\end{prop}
\begin{dem}
Puisque $A$ et $A'$ engendrent le m\^eme sous-groupe du groupe de Brauer de $K$, on sait (cf. introduction) que $SB(A)\times \P^{n-1}$ est birationnelle \`a $SB(A')\times \P^{n-1}$.  D'apr\`es le corollaire pr\'ec\'edent, $SB(r,A)$ (resp. $SB(r,A')$) est birationnelle \`a $SB(A) \times \P^s$ (resp. $SB(A') \times \P^s$) avec $s=r(n-r)-n+1$. Puisqu'on a $s \geq n-1$ au vu des hypoth\`eses faites sur $n$ et $r$ (v\'erification laiss\'ee au lecteur), on en d\'eduit que $SB(A) \times \P^s$ et $SB(A') \times \P^s$ sont birationnellement isomorphes, cqfd.
\end{dem}

\signmf

\end{document}